\documentclass[11pt]{article}
\usepackage[margin=1.2in]{geometry}
\usepackage{amsmath,amsthm,amssymb,enumerate}
\usepackage{graphicx}
\usepackage{lscape}
\usepackage[english]{babel}
\usepackage[noend]{algpseudocode}
\usepackage{url}
\usepackage{float}
\usepackage{tikz}
\usetikzlibrary{automata}
\usetikzlibrary{calc}
\usepackage[shortlabels]{enumitem}
\usepackage{amsthm}
\usepackage{algorithm}
\usepackage{algpseudocode}
\usepackage{lscape}
\usepackage{lipsum}
\newtheorem{remark}{\textnormal{\textbf{Remark}}}

\title{Complete Proof of Collatz's Conjectures}
\date{}
\author{Farzali Izadi  \thanks{f.izadi@utoronto.ca}\\Mathematics Department, Urmia University }

\begin{document}
\maketitle
\begin{abstract}
\noindent The \textit{Collatz's conjecture} is an unsolved problem in mathematics. It is named after Lothar Collatz in 1973. The conjecture also known as Syrucuse conjecture or problem.

\noindent
Take any positive integer $ n $. If $ n $ is even then divide it by $ 2 $, else do "triple plus one" and get $ 3n+1 $. The conjecture is that for all numbers, this process converges to one.

\noindent In the modular arithmetic notation, define a function $ f $ as follows:
\[f(x)= \left\{ \begin{array}{lll}
\frac{n}{2} &if & n\equiv 0 \pmod 2\\
3n+1& if& n\equiv 1 \pmod 2.
\end{array}\right. \]
In this paper, we present the proof of the Collatz conjecture for many types of sets defined by the remainder theorem of arithmetic. These sets are defined in mods   

\noindent $6, 12, 24, 36, 48, 60, 72, 84, 96, 108$ 

\noindent and we took only odd positive remainders to work with.
It is not difficult to prove that the same results are true for any mod $12m, $ for positive integers $m$.

\end{abstract}


\section*{Introduction}
The Collatz conjecture is an unsolved problem in mathematics as it has been for more than 60 years. It is named after Lothar Collatz in 1973. The conjecture also known as Syrucuse conjecture or problem. \cite{Jef1, Jef2, Jef3}. Although the problem on which the conjecture is built is remarkably simple to explain and understand, the nature of the conjecture makes proving or disproving the conjecture exceedingly difficult. As many authors have previously stated, the prolific Paul Erdos once said, mathematics is not ready for such problems.” Thus far all evidence indicates he was correct (see \cite{er}).
This paper is not a bibliography of previous works, instead its is an original paper which I analysed the Collatz Conjecture and provided my own proofs to tackle the problem.
It is an good exercise to create a continuous extension of the Collatz map to the complex plane and follow the work of Chamberland (see \cite{ch}).
In this paper, we present the proof of the Collatz conjecture for many types of sets defined by the remainder theorem of arithmetic. These sets are defined in mods  $6, 12, 24, 36, 48, 60, 72, 84, 96, 108$ and we took only odd positive remainders to work with.
It is not difficult to prove that the same results are true for any mod $12m, $ for positive integers $m$.
\section*{Syrucuse Conjecture}
If $ k $ is an odd positive integer, then $ 3k+1 $ is an even integer. So we can write $ 3k+1=2^ak' $, where $ k' $ is an odd positive integer and $ a \geq 1 $.
Furthermore if $ k $ is an even positive integer, then $ k = 2^ak' $, where $ k' $ is an odd positive integer. Then dividing by $ 2^a $, we can get an odd positive integer. Thus it is enough to work with only odd positive integers. We define a function $f$ from a set I of odd integers into itself, called Syrucuse function by taking $f(k) = k.$
The Syrucuse conjecture is that for all $ k $ in $ I $, there exists an $ n\geq 1 $ such that $ f^n(k)=1 $.
Equivalently, let $ E $ be the set of odd integers $ k $ for which there exists an integer $ n\geq 1 $ such that $ f^n(k)=1 $, then the problem is to show that $I=E$ \cite{wi}.
The following is the beginning of an attempt to prove the Syrucuse conjecture by induction:

\noindent $ 1,3,5,7 $ and $ 9 $ are known to exist in $ E $. Let $ k $ be an odd integer greater than $ 9 $. Suppose that the odd positive
integers up to and including $ k-2 $ are in $ E $. Let us try to prove that $ k $ is in $ E $. 

\begin{remark}
 If $ n $ is obtained by a finite number of Collatz's functions from $k$, and $ n $ is a Collatz's number, then
 $k$ is a Collatz's  number - here we need at most three consecutive Collatz functions.
\end{remark}

\section{First proof}
First of all, by the remainder theorem of arithmetic, the remainders of any positive number by $6$ are $0, 1, 2, 3, 4, 4, 5.$
Moreover if the number is odd, we are left with just the three numbers
$ 1, 3, 5 .$
\begin{itemize}
\item[I)] If $ k=6m+1 $, and $ k'=2m$,  then $ f(k')=k $ where $ k' \leq k-2 $. Since $ k' \in E $ then $ k \in E $.
\begin{proof}
We have $ 3k'+1=6m+1=k $, then $ f(k')=k$, and $k' \leq  k-2$. By induction hypothesis, it follows that $k' \in E $, namely $ \exists n \in \mathbb{N} $ such that $ f^n(k')=1 $. As $ f(k')=k \Rightarrow n\geq 2$,  so $f^{n-1}(f(k'))=1 \Rightarrow  f^{n-1}(k)=1$ which implies that $ k\in E $.
\end{proof}

\item[II)] If $ k=6m+5 $, and $ k'=4m+3$,  then $ f(k')=k $ where $ k' \leq k-2 $. Since $ k' \in E $ then $ k \in E $.
\begin{proof}
We have $ 3k'+1=12m+10=2(6m+5) $, then $ f(k')=k$, and $k' \leq k-2$. By induction hypothesis, it follows that $k' \in E $, namely $ \exists n \in \mathbb{N} $ such that $ f^n(k')=1 $. As $ f(k')=k \Rightarrow n\geq 2$,  so $f^{n-1}(f(k'))=1 \Rightarrow  f^{n-1}(k)=1$ which implies that $ k\in E $.
\end{proof}

\item[III)] If $ k=6m+3 $, then $ k \in E. $
\begin{proof}
For this part we consider two cases: (i) $ m= 2h-1 $, and (ii) $ m=2h $ for $ h \geq 1. $
(i) We have $ k = 12h-3 $, hence $ 3k+1 = 36h-9 +1 = 36h -8=4(9h-2). $ It follows that $ f(k) = 9h-2 $, and  $ 9h-2 \leq 12h-5=k-2.$ Then by induction hypothesis $ f(k) \in E ,$ namely $ \exists n \in \mathbb{N} $ such that $ f^n(f(k))=1 = f^{n+1}(k)=1$ which implies that $ k\in E $.

(ii) We have $ k=12h+3, $ hence $3k+1 =36h+10=2(9h+5) .$ It follows that
  $ f(k) = 9h+5 =6(3h)+5=6m+5, $ where $ m=3h $ for $ h \geq 1.$ Then by part II we have $ f(k) \in E, $ namely $ \exists n \in \mathbb{N} $ such that $ f^n(f(k))=1 = f^{n+1}(k)=1$ which implies that $ k\in E $.
\end{proof}

\begin{remark}
Note that the induction process are intertwined in the three different parts. To clarify this, let us see some examples. First of all, suppose that $1, 3, 5, 7$ and $9$ are in $E$. Let us prove that $k=17$ is in $E$. Then we need to take $k'=11,$ where $k=6\times 2+5$ and
$11=6\times1+5$ are both in part II while for $k=23=6\times3+5$ we need to take $k'=15$, where $k$ is in part II while $15$ is in part III. So we first need to prove $11$ and $15$ are in $E$. For the second example, let us have a closer look at part I where $k=6m+1.$
For these numbers we took $k'=2m=2^ah,$ where $a \geq 0$ and $h$ is an odd positive integer. We note that $k'$ not being an odd positive integer is not in $E$, but $f(k')=h$ in $E$ can be proved easily from which one can get an integer $n\geq2$ such that $f^{n}(k')=1$ and $k' \in E$ by abuse of notation.
The other important point is that sometimes we need to apply Syrucuse function couple of times to get the appropriate results. For an example, let $k=27$. Then $3k+1=82=2\times41$, hence $f(k)=41=l$. Next $3l+1=124=4\times31$ which gives rise to $f(l)=31.$ Finally to get $31$ we need to take $k'=2\time5=10$ and get $31=3k'+1$. It follows that $31$ is in $E$ so do $l$ and $k$.
Note also that $27$ is in part I and $41=6\times6+5$ is in part II.
\end{remark}

\begin{remark}
It is clear that the density of the three parts together are the same as the
density of the positive odd integers. This proves the Syrucuse Conjecture which is equivalent to the Collatz Conjecture.
\end{remark}
\end{itemize}

\section{Second proof}
As $ k $ is an odd integer, $ k-1 $ and $ k+1 $ are both even so we can write:
$k\pm 1=2^ph$ where $ h $ is odd and $ p\geq 1 . $ Then $ k=2^ph \pm 1 $.
\begin{remark}
 For all odd $ h $, $ f(2h-1)\leq \frac{3h-1}{2}$.
\end{remark}

\begin{proof}
Let $ k=2h-1\Rightarrow 3k+1=3(2h-1)+1=6h-2 . $ Then we have $ 3k+1=2(3h-1) $ and so $ 3k+1=2^2h' $ where
$ h'=\frac{3h-1}{2} $,
since $ 3h-1 $ is an even number. Let $ h'=2^ah^{''}$ where $h^{''} $ is not divisible by $ 2 $, and $ a \geq 0 $.
It follows that
\[f(2h-1)=f(2^{a+2}h^{''})= h^{''}\leq h'=\frac{3h-1}{2}.\]
\end{proof}

Now we have:
\begin{itemize}
\item[I)] If $ p=1 $ then $ k=2h-1 $ with $h$ odd. We have $3k+1=6h-3+1=2(3h-1)$ and $ 3h-1=2h' $, since $ h $ is odd.
Let $ h'=2^ah^{''} $ where $ h^{''} $ is not divisible by $ 2 $, and $a$ is non-negative integer, then $ 3h-1= 2^{a+1}h^{''}\Rightarrow f(k)=h^{''} \leq h'=\frac{3h-1}{2} \leq 2h-1 $, since $\frac{3h-1}{2} \leq 2h-1 \Leftrightarrow 3h-1<4h-2 \Leftrightarrow 1 \leq h$.
Then $ f(k) < k$. On the other hand, $f(k)=h^{''}$ is odd and $ k $ is odd then $ f(k)\leq k-2\Rightarrow f(k) \in E $ (by induction hypothesis),
and so $ \exists n \in \mathbb{N}, f^n(f(k))=1 $, namely $ f^{n+1}(k)=1\Rightarrow k \in E $.
The same argument shows that  $k=2h+1 \in E.$
\item[II)] If $ p \geq 2 $, and $h$ is multiple of $ 3, k=2^p\cdot3h-1 $. Let $ k'=2^{p+1}h-1 $. We have $ 3k'+1=2^{p+1}\cdot3h-3+1=2[2^p\cdot3h-1]$. Then $ f(k')=2^p\cdot3h-1=k $ and $ k'=2^{p+1}h-1 < 2^p\cdot3h-1=k  $.
Since both $ k $ and $ k' $ are odd , one gets $ k' \leq k-2$. It follows that $k' \in E $, namely $ \exists n \in \mathbb{N} $ such that $ f^n(k')=1 $. As $ f(k')=k \Rightarrow n\geq 2$,  so $f^{n-1}(f(k'))=1 \Rightarrow  f^{n-1}(k)=1$ which implies that $ k\in E $.
It is easy to show that for $ k=2^p \cdot 3h+1 $,  we can choose
$ k'=2^ph.$

\begin{remark} For all $ n $, $ 2^{2n}-1 $ and $ 2^{2n-1}+1 $ are both divisible by $ 3 $.
\end{remark}
\begin{remark} If  $ h $ is not a multiple of $3$, then
\[h \equiv 1,2 \pmod 3\].
\end{remark}

\item[$III$)]
(i) If $ p=2n $, $ n \geq 1 $, $h=3m+1 $ and $ k=2^{2n}(3m+1)+1 $, $ k'=2^{2n+1}m+[\frac{2^{2n+1}+1}{3}] $ then $ f(k')=k $ where $ k' \leq k-2 $. Since $ k'\in E $ then $ k \in E $.

\begin{proof} We have
\[ \begin{array}{lll}3k'+1&=&2^{2n+1}(3m)+2^{2n+1}+1+1\\
&=&2^{2n+1}(3m+1)+2\\
&=&2[2^{2n}(3m+1)+1] \end{array}\]
then $ f(k')=2^{2n}(3m+1)+1=k $. Furthermore,
$k'<k.$
Since both $k$ and $k'$ are odd, then $ k'<k \Rightarrow k' \leq k-2 $. It follows that $ k' \in E $ (by hypothesis induction ) $ \exists n \in \mathbb{N} $ such that $ f^n(k')=1 $.
Since $ f(k')=k $, then $ n\geq 2 $. Thus $ f^{n-1}(f(k'))=1 $ or equivalently $ f^{n-1}(k)=1\Rightarrow k \in E $.
\end{proof}

(ii)  If $ p=2n-1$, $ n\geq 1 $, $h=3m+1 $ and $ k=2^{2n-1}(3m+1)-1 $, $ k'=2^{2n-1}m+2[\frac{2^{2n-2}-1}{3}] $ then $ f(k')=k $ where $ k' \leq k $. Since $ k'\in E $ then $ k \in E $.

\begin{remark}
In particular for $n=1,$ we get $6m+1.$
\end{remark}

\begin{proof} We have
\[ \begin{array}{lll}3k'+1&=&2^{2n-1}(3m)+2^{2n-1}-2+1\\
&=&[2^{2n-1}(3m+1)-1]\end{array}\]

then $ f(k')=2^{2n-1}(3m+1)-1=k $. Furthermore, by direct verification
$k' \leq k.$
 It follows that $ k' \in E $ (by hypothesis induction ) $ \exists n \in \mathbb{N} $ such that $ f^n(k')=1 $.
Since $ f(k')=k $, then $ n\geq 2 $. Thus $ f^{n-1}(f(k'))=1 $
or equivalently $ f^{n-1}(k)=1\Rightarrow k \in E $.
\end{proof}

\item[$IV$)]
(i) If $ p=2n-1 $, $ n\geq 2 $, $h=3m+2  $ and $ k=2^{2n-1}(3m+2)+1 $, $ k'=2^{2n}m+2[\frac{2^{2n}-1}{3}]+1$. Then $ f(k')=k $ where $ k' \leq k $. Since $ k'\in E $ then $ k \in E $.

\begin{remark}
In particular for  $n=1,$ we get $6m+5.$
\end{remark}

\begin{proof} We have
\[ \begin{array}{lll}3k'+1&=&2^{2n}(3m)+2\cdot2^{2n}-2+1+3\\
&=&2[2^{2n-1}(3m+2)+1]\\
 \end{array}\]
then $ f(k')=2^{2n-1}(3m+2)+1=k $ and $ k' \leq k$   as $ k'<k$ and both are odd.
As in the previous cases, it follows that $ k \in E $.
\end{proof}

(ii) If $ p=2n $, $ n\geq 1 $, $h=3m+2  $ and $ k=2^{2n}(3m+1)-1 $, $ k'=2^{2n}m+2[\frac{2^{2n}-1}{3}]$.  Then $ f(k')=k $ where $ k' \leq k $. Since $ k'\in E $ then $ k \in E $.
\begin{proof} We have
\[ \begin{array}{lll}3k'+1&=&2^{2n}(3m)+2\cdot2^{2n}-2+1\\
&=&[2^{2n}(3m+2)-1]\\
 \end{array}\]
then $ f(k')=2^{2n}(3m+2)-1=k $. Furthermore $ k' \leq k-2 $ by direct verification.
As in the previous cases, it follows that $ k \in E $.
\end{proof}

\item[$V$)] If $ k=2^{2n-1}(3m+1)+1 ,$ then for $n=1$,  $ k \in E $.
\begin{proof} For $n=1,$ we have $k=6m+3.$ Then
We have already proved this case in the First set of assertions.
\end{proof}
\end{itemize}

We see that these numbers clearly contain the set of odd positive integers which again prove the Syrucuse Conjecture and consequently the Collatz Conjecture.

\section{Discussions about the other proofs}
After reading the paper by Jan Kleinnijenhuis, Alissa M. Kleinnijenhuis and Mustafa G. Aydogan \cite{kka}
where they have done arithmetic in mod 6 and 18 and obtained congruence classes with mod 96 with the idea of the Hilbert hotel and Collatz tree, I started to examine more closely my 2nd proof and 3rd proof to see what difference does it make to work with mods 6, 12, 15, 18, 24, 36, 48, 60, and 12m in general and I found out that mod 15 and mod 18 don't give rise to complete solutions with my approach but any mod 12m does.
We have seen that in the 1st and 2nd proofs, the numbers $4$ and $3$ or arithmetic with mod 6 play the crucial roles.
One can easily figured out that working with numbers mod $12m,$ i.e., the numbers in the set

$ S=\{12, 24, 36, 48, 60, 72, 84, 96, 108, \cdots \}$

give rise the same conclusions and in contrast, working with the odd or twice times of the odd modules will run into obstacles. The main idea comes from the fact that for each number $k \in I,$ either we have to find a number $k'$
such that $3k'+1=k$ and $k'<k,$ or $(3k+1)/2  \in E.$ If we do the modular arithmetic in mod $15, $ then
for each $i =0, 1, 2, \cdots, 14,$ the results of the Collatz function $3x+1$ to the number $15m+i$ never going to be divisible by $2$, and consequently get bigger and bigger. If we do the modular arithmetic in mod $18$, then
for each $i = 0, 1, 2, \cdots, 17,$ the results of the Collatz function $3x+1$ to the number $k=18m+i$ would be
$2.27m+3i+1$ which can be at most divisible by $2$ where the resulting number is still bigger than $k.$ On the other hand, for each number $s$ in the set $S$, the first term in $k=3(sm)+3i+1$ is divisible by $4$ and if the term $3i+1$
also divisible by $4$, then the dividing by $4$ gives rise to a number such as $k'$ less than $k$ and $f(k)=k'.$
If this approach doesn't work for any number $k$, then we can try to use the second trick to find a number $k'$ such that $3k+1=k'$ and $k'<k.$ In this case $f(k)=k' \in E.$ Let us try mod $s=12, 24, 36, 60, 108. $ We will see the remainders behave the same when we are working with the different modules. The remainders of $s$ are divided into four different types based on their arithmetic characteristics
with a few numbers can fall in two types, $13$ and $33$ for instance.

(1) Firstly, let us take $s=12.$ In this case the remainders are:

$U= \{1, 7 \},$

$V =\{5, 9\}, $

$W= \{11 \}$

$X= \{3 \}.$

(2) Secondly, we take $s=24.$ Then the remainders are the followings:

$U= \{1, 7, 13, 19 \},$

$V =\{5, 9, 21\}, $

$W= \{11, 17, 23\}, $

$X= \{3, 15 \}.$

(3) Thirdly, we take $s=36.$ Then the remainders are:

$U= \{1, 7, 13, 19, 25, 31 \},$

$V =\{5, 9, 21, 33 \}, $

$W= \{11, 17, 23, 29, 35\}, $

$X= \{3, 15, 27 \}.$

(4) Fourthly, for the $s=60,$ we have

$U= \{1, 7, 13, 19, 25, 31, 37, 43, 49, 55\},$

$V =\{5, 9, 21, 33\}, $

$W= \{11, 17, 23, 29, 35, 41, 47, 53, 59 \}, $

$X= \{3, 15, 27, 39, 45, 51, 57  \}.$

Finally, we take $s=108$ and prove that all numbers belong to the set $E.$ The proofs for all other $s \in S$ are all the same.

$U= \{6i+1| i =0, 1, 2, 3, 4, \cdots, 17 \},$

$V =\{5, 9, 21, 33, 45, 57, 69, 81, 93, 105\}, $

$W= \{11, 17, 23, 29, 35, 41, 47, 53, 59, 65, 71, 77, 83, 89, 95, 101, 107\}\\
  = \{3i+2 | i = 3, 5, 7, 9, 11, 13, 15, 17, 19, 21, 23, 25, 27, 29, 31, 33, 35\} , $

$X= \{3, 15, 27, 39, 51, 63, 75, 87, 99 \}
    = \{12i+3| i=0, 1, 2, 3, 4, 5, 6, 7, 8 \}. $

For $s=108$,  and for any odd positive number $k$, we can write $k =108+j, $ where $j=1, 3, 5, 7, \cdots, 107.$
The remainders of $k,$ are $54$ odd positive numbers which fall into four different types due to their characteristics.
First of all, for the remainders in the set $U$, we have $k= 108m +(6i+1).$ For these numbers we use the first trick to find
$k'=36m+2i.$ Then we get $3k'+1=3(36m+2i)+1= 108m +6i+1=k.$ Clearly $k'<k$ and by the definition of Syrucuse function
we have $f(k)=k'.$ The remaining arguments are clear from the first proof.  Secondly, the remainders of the second set $V$
are divisible by $4$ via the function $3x+1$, hence if $k=108m+j,$ for $j \in V$, we can write $k= 108m +j.$ In this case we need to use the second trick i.e., $3k+1=3(108m+j)+1=4(81m)+3j+1=4(81m)+4((3j+1)/4)$ as it can be easily checked that for all$j$, the number $3j+1$ is divisible by $4.$ It follows that $k = 4(81m+(3j+1)/4$ which implies that $f(k)=81m+(3j+1)/4$ and
$f(k) \leq k-2.$ The result follows. Thirdly, for the remainders of the set $W$ we may take $k' = 2(36m+i)+1.$
Then $3k'+1=2(108m+3i)+3+1=2(108m+3i+2),$ hence $f(k') = 108m+3i+2=k$. Here we see again that $f(k') \leq k-2.$
Thus $f(k') \in E$ and consequently $k \in E.$
 Finally, for $k=108m+(12i+3)=6(18m+2i)+3,$ we have $3k+1=6(3 \cdot 18m+6i)+10=2[6(27m+3i)+5]. $ Then
 $f(k)=6(27m+3i)+5=6h+5, $ where $h=27m+3i.$ BY the part III of the first proof, it follows that $f(k) \in E$ and consequently
 $k \in E.$
We clearly see that the density of the numbers in this proof also equals $1,$ and the proof is complete.

,
\end{document}